\documentclass[12pt]{article}

\usepackage{amsthm,amsmath,amssymb}

\usepackage{graphicx}

\usepackage[colorlinks=true,citecolor=black,linkcolor=black,urlcolor=blue]{hyperref}

\newcommand{\doi}[1]{\href{http://dx.doi.org/#1}{\texttt{doi:#1}}}

\theoremstyle{plain}
\newtheorem{theorem}{Theorem}

\newtheorem{corollary}[theorem]{Corollary}
\newtheorem{proposition}[theorem]{Proposition}

\theoremstyle{definition}
\newtheorem{definition}[theorem]{Definition}

\newtheorem{conjecture}[theorem]{Conjecture}

\theoremstyle{remark}

\title{\bf On-Line Choice Number of Complete Multipartite Graphs: an Algorithmic Approach}

\author{Fei-Huang Chang\thanks{ROC National Science Council grant NSC101-2115-M-003-007.}\\
\small Division of Preparatory Programs\\[-0.8ex]
\small for Overseas Chinese Students\\[-0.8ex]
\small National Taiwan Normal University\\[-0.8ex]
\small New Taipei, Taiwan\\
\small\tt feihuang0228@gmail.com\\
\and
Hong-Bin Chen\\
\small Institute of Mathematics\\[-0.8ex]
\small Academia Sinica\\[-0.8ex]
\small Taipei, Taiwan\\
\small\tt andanchen@gmail.com\\
\and
Jun-Yi Guo\thanks{ROC Ministry of Science and Technology grant 103-2511-S-003-006-.}\\
\small Department of Mathematics\\[-0.8ex]
\small National Taiwan Normal University\\[-0.8ex]
\small Taipei, Taiwan\\
\small\tt junyiguo@gmail.com\\
\and
Yu-Pei Huang\thanks{ROC Ministry of Science and Technology grant 103-2811-M-214-001 and 103-2632-M-214-001-MY3-1.}\\
\small Department of Applied Mathematics\\[-0.8ex]
\small I-Shou University\\[-0.8ex]
\small Kaohsiung, Taiwan\\
\small\tt yphuang@isu.edu.tw
}

\begin{document}

\maketitle


\begin{abstract}
  This paper studies the on-line choice number of complete multipartite graphs with independence
number $m$. We give a unified strategy for every prescribed $m$. Our main result leads to several
interesting consequences comparable to known results. (1) If $\displaystyle
k_1-\sum_{p=2}^m\left(\frac{p^2}{2}-\frac{3p}{2}+1\right)k_p\geq 0$, where $k_p$ denotes the number
of parts of cardinality $p$, then $G$ is on-line chromatic-choosable. (2) If $
|V(G)|\leq\frac{m^2-m+2}{m^2-3m+4}\chi(G)$, then $G$ is on-line chromatic-choosable. (3) The
on-line choice number of regular complete multipartite graphs $K_{m\star k}$ is at most
$\left(m+\frac{1}{2}-\sqrt{2m-2}\right)k$ for $m\geq 3$.

\end{abstract}

\section{Introduction}

On-line list coloring was introduced by Schauz \cite{Schauz09} in the context of the Paint-Correct
game played on a graph. Here we adapt the on-line list coloring version introduced by Zhu
\cite{Zhu09} as follows.

\begin{definition}
Given a graph $G$ and an integer-valued function $f$ on $V(G)$, the on-line $f$-list coloring of
$G$ is a two-players game, say Alice and Bob, played on $G$. In the very beginning, all vertices
are uncolored. In the $i$th round, Alice marks a nonempty subset $V_i$ of remaining uncolored
vertices and assigns color $i$ as a permissible color to each vertex of $V_i$. Then Bob chooses an
independent set $X_i$ contained in $V_i$ and colors all vertices of $X_i$ the color $i$. The game
goes round by round. If at the end of some round there is a vertex $v$ which has been assigned
$f(v)$ permissible colors, i.e., has been marked $f(v)$ times, but is not yet colored by Bob, then
Alice wins the game. Otherwise, Bob wins, i.e., in the end each vertex $v$ is colored by Bob before
running out of $f(v)$ permissible colors.
\end{definition}

Given an integer-valued function $f$ defined on $V(G)$, we say that $G$ is \emph{on-line
$f$-choosable} if Bob has a winning strategy for the on-line $f$-list coloring game on $G$ no
matter how Alice plays; particularly, if $f(v)$ is a constant $k$ for all $v\in V(G)$, then we say that
$G$ is \emph{on-line $k$-choosable}. Denoted by $\chi_p(G)$, the \emph{on-line choice number} of $G$
is the minimum number $k$ such that $G$ is on-line $k$-choosable.

The conventional list coloring, introduced by Vizing \cite{vizing} and independently by Erd\H{o}s,
Rubin and Taylor \cite{ERT80}, is a special case that Alice shows Bob the full lists in the very
beginning of the on-line list coloring game. So Bob has a winning strategy for the list coloring if
he has one for the on-line list coloring game. Let $\chi(G)$ and $\chi_\ell(G)$ denote the \emph{
chromatic number} and \emph{choice number} of a graph $G$, respectively. In general, we have
$\chi(G)\leq\chi_\ell(G)\leq\chi_p(G)$ for any $G$.

It is known that $\chi_\ell(G)-\chi(G)$ can be arbitrarily large; see \cite{GS} for an example that
demonstrates complete bipartite graphs $G$ having $\chi_\ell (G)$ arbitrarily large but
$\chi(G)=2$. An interesting question is whether $\chi_p(G)-\chi_\ell(G)$ can be arbitrarily
large.  To the best of our knowledge, the problem is still open. Although there exist a few graphs
$G$ with $\chi_p(G)> \chi_\ell(G)$, the largest gap known up to now is $1$ as shown in \cite{Zhu09}.

Another interesting question raises naturally: for which graphs $G$ does $\chi(G)=\chi_\ell(G)\\ =\chi_p(G)$? There
have been many studies on graphs satisfying $\chi_\ell(G)=\chi(G)$ (see \cite{galvin,JT,KW01,ohba}
and references therein); such a graph is called \emph{chromatic-choosable}. Likewise, a graph $G$ is
called \emph{on-line chromatic-choosable} if $\chi_p(G)=\chi(G)$. Ohba \cite{ohba} conjectured that
for all graphs $G$ with $|V(G)|\leq 2\chi(G)+1$, $G$ is chromatic-choosable; recently, this has
been proved by Noel, Reed and Wu \cite{RNW13+}. Let $K_{n_1\star k_1,n_2\star k_2,\ldots,n_s\star
k_s}$ denote the complete multipartite graph where $k_i$ partite sets are of size $n_i$ for
$i=1,2,\ldots,s$. For short, we shall simplify $n_i\star 1$ as $n_i$ (for example $K_{3\star
2,4}=K_{3\star 2,4\star 1}$). In view of the fact that $K_{2\star k,3}$ is not on-line
chromatic-choosable for $k\geq 2$ \cite{KKLZ12}, Huang, Wong and Zhu \cite{HWZ11} slightly modified
the Ohba's conjecture to its on-line version.

\begin{conjecture}\label{online}\cite{HWZ11}
Every graph $G$ with $|V(G)|\leq 2\chi(G)$ is on-line chromatic-\\
choosable.
\end{conjecture}
We remark that to prove the on-line Ohba's conjecture, it suffices to prove it for complete
$\chi(G)$-partite graphs $G$ since adding edges does not reduce the on-line choice number. Conjecture \ref{online} has been verified for complete multipartite graphs with a small independence number.
Using the Combinatorial Nullstellensatz, the authors \cite{HWZ11} proved that $K_{2\star k}$ is
on-line chromatic-choosable. Recently, Kim et al \cite{KKLZ12} gave an algorithmic proof for $K_{2\star
k}$ and later Kozik, Micek and Zhu \cite{KMZ13+} extended the case to complete multipartite graphs with
independence number at most 3.

This paper focuses on complete multipartite graphs with independence number $m$.  In Section 2, we
generalize the algorithmic methods in \cite{KKLZ12,KMZ13+} and give a unified strategy for the
on-line choice number of graphs with any prescribed $m$. Our main result provides a sufficient
condition on $f$ for graphs being on-line $f$-choosable by partitioning vertices into independent sets in a systematic
way. It is a broadly applicable tool which leads to several interesting
consequences comparable to known results. Section 3 presents some immediate consequences.


%
%
%


\section{Main Result}

This section starts with some notations and definitions. Throughout the rest of this paper, we shall use ``part'' instead of ``partite set'' for short, and let $m$ be a fixed positive integer. Consider a complete multipartite graph $G$
with part size at most $m$, i.e. independent number at most $m$. Let $\Pi=\{\underline{X}_{m-1}$,
$\underline{X}_{m-2}, \ldots, \underline{X}_{1}$, $\overline{X}_{2}, \ldots, \overline{X}_{m}\}$ be a partition of parts of
$G$, where $\overline{X}_{p}$ is a family consisting of parts of size exact $p$ for $2\leq p\leq m$ and
$\underline{X}_{p}$ a family consisting of parts of size at most $p$ for $1\leq p\leq m-1$. Particularly,
$\underline{X}_{1}$ contains only parts of size $1$.  Notably, the partition is not unique since a part
of size $p$ can belong to either $\overline{X}_{p}$ or $\underline{X}_{j}$ for some $j\leq p$. Let $u_p$ and
$\ell_p$ denote the number of parts of $\overline{X}_p$ and $\underline{X}_p$, respectively, i.e.,
$|\overline{X}_p|=u_p$ and $|\underline{X}_p|=\ell_p$. For each family $\underline{X}_{p}$, we use the second
coordinates to denote parts in
            the family, e.g.,
            $\underline{X}_{p}=(\underline{X}_{p,1},\ldots,\underline{X}_{p,\ell_p})$ where
            $\underline{X}_{p,i}$ means the $i$th part in $\underline{X}_{p}$.   When it comes to vertices in a family of parts, we shall use the notation
$V(\cdot)$ to avoid confusion. Given a function $f: V(G)\rightarrow \mathbb{N}$, for
$y=1,2,\ldots,m,$ define
$$F(y)\equiv \min\left\{ \sum_{v\in V(Y)}f(v): Y\subseteq X\in \bigcup_{p=2}^m\overline{X}_p \mbox{ and }
|Y|=y\right\}.$$ Notice that $\displaystyle F(1)=\min_{v\in \bigcup_{p=2}^mV(\overline{X}_p)} f(v).$ For $j=1,2,\ldots, m$,
define
$$S(j)\equiv\sum_{p=2}^ju_p+\sum_{p=j+1}^m(p-j)u_p.$$
Let $\displaystyle\alpha(1)=\sum_{p=2}^mu_p$ and $\beta(1)=0$, and define recursively that for $j=2,\ldots,m$,
$$\alpha(j)\equiv \alpha(j-1)+S(j-1)  \mbox{ and }$$
$$\beta(j)\equiv\beta(j-1)+\sum_{p=1}^{m-j+1}|V(\underline{X}_p)|.$$
Then it can be easily expressed as
\begin{equation}\label{eq1}\alpha(j)=\sum_{p=2}^j\left(j+\frac{p^2}{2}-\frac{3p}{2}+1\right)u_p+\sum_{p=j+1}^m\left(\frac{j}{2}-\frac{j^2}{2}+pj-p+1\right)u_p
\end{equation} and
\begin{equation}\label{eq2}\beta(j)=\sum_{p=1}^{m-j}(j-1)|V(\underline{X}_p)|+\sum_{p=m-j+1}^{m-1}(m-p)|V(\underline{X}_p)|\end{equation}
for $j=1,\ldots,m$.

The following propositions are elementary but useful observations.

\begin{proposition} For $j=1,\ldots,m-1$, let $\displaystyle S(m-j)=\sum_{p=2}^m s_pu_p$. Then we have $1=s_2= s_3=\cdots=
s_{m-j+1}<s_{m-j+2}<\cdots< s_{m-1}<s_m$. \label{pro1} \end{proposition} \proof The proof follows immediately
by definition.\qed

\begin{proposition} For $j\geq 2$, let $\displaystyle\alpha(j)+\beta(j)=\sum_{p=2}^m
a_pu_p+\sum_{p=1}^{m-1}b_p|V(\underline{X}_p)|$ and $a_1=b_1$. The following is true.\begin{enumerate}
\item[(i)] $b_1=b_2=\cdots=b_{m-j+1}=j-1$.
\item[(ii)] For any integers $s$ and $t$ with $1\leq s<t\leq m$, we have $a_{t}-a_s\geq
t-s$. In particular, if $t>j$, then $a_{t}-a_s\geq j-1$ and if $t\leq j$, then $a_{t}-a_s\geq t-2$.
\item[(iii)] For $p=2,3,\ldots,m$, we have $a_p\geq\max\{j,p\}$.  
\end{enumerate}
\label{pro2}
\end{proposition}

\proof The proof of (i) is trivial.
We first prove (ii) and consider the following cases.
\begin{enumerate}
  \item[$s=1$]
        \begin{enumerate}
          \item If $2\leq j<t\leq m$, then
                        $(a_t-a_s)-(t-s)=(a_t-a_1)-(t-1) \\
                                        =\left(\left(\frac{j}{2}-\frac{j^2}{2}+tj-t+1\right)-(j-1)\right)-(t-1) \\
                                        =\frac{j}{2}-\frac{j^2}{2}+t(j-2)+3\geq
                                        \frac{j}{2}-\frac{j^2}{2}+(j+1)(j-2)+3 \\
                                        =\frac{j(j-1)+2}{2}\geq 0$. Thus, $a_t-a_s\geq
                 \max\{t-s, j-1\}$.
          \item If $1<t\leq j\leq m$, then $(a_t-a_s)-(t-s)=(a_t-a_1)-(t-1)=((j+\frac{t^2}{2}-\frac{3t}{2}+1)-(j-1))-(t-1)=\frac 12{(t-2)(t-3)}\geq 0$
          (as $t\in \mathbb{N}$). Therefore, $a_t-a_s\geq \max\{t-s, t-2\}$.
        \end{enumerate}
  \item[$s>1$]
        \begin{enumerate}
          \item If $2\leq j<s<t\leq m$, then we have $a_t-a_s=(\frac{j}{2}-\frac{j^2}{2}+tj-t+1)-(\frac{j}{2}-\frac{j^2}{2}+sj-s+1)=(t-s)(j-1)\geq \max\{t-s,
          j-1\}$.
          \item If $2\leq s\leq j<t\leq m$, then
                    \begin{align*} (a_t-a_s)-(t-s) &=\textstyle((\frac{j}{2}-\frac{j^2}{2}+tj-t+1)-(j+\frac{s^2}{2}-\frac{3s}{2}+1))-(t-s) \\
                                                   &=\textstyle\frac 12(2t(j-2)-j-j^2-s^2+5s)\\
                                                   &\geq \textstyle\frac 12(2(j+1)(j-2)-j-j^2-s^2+5s) \\
                                                   &=\textstyle\frac 12{((j+s-3)(j-s)+2(s-2))}\geq 0.
                    \end{align*}
                Besides,
                    \begin{align*}(a_t-a_s)-(j-1) &\textstyle=\left(\left(\frac{j}{2}-\frac{j^2}{2}+tj-t+1\right)-\left(j+\frac{s^2}{2}-\frac{3s}{2}+1\right)\right)-(j-1) \\
                                                  &=\textstyle\frac 12 \left(2t(j-1)-j^2-3j+2-(s^2-3s)\right)\\
                                                  &\geq \textstyle\frac 12\left(2(j+1)(j-1)-j^2-3j+2-(s^2-3s)\right) \\
                                                  &=\textstyle\frac
                12\left((j^2-3j)-(s^2-3s)\right)\geq 0,\end{align*} where the last
                inequality holds for $j\geq s\geq 2$. Therefore, $a_t-a_s\geq
                \max\{t-s,j-1\}$.
          \item If $2\leq s<t\leq j\leq m$, then $a_t-a_s=(j+\frac{t^2}{2}-\frac{3t}{2}+1)-(j+\frac{s^2}{2}-\frac{3s}{2}+1)=\frac 12(t-s)(t+s-3)$.
          If $s=t-1$, then $a_t-a_s=t-2\geq t-s$. If $s<t-1$, then $t-s\geq 2$ and
          $a_t-a_s\geq t+s-3\geq t-1>t-2\geq t-s$. In either case, $a_t-a_s\geq \max\{t-s,
          t-2\}$, as desired.
        \end{enumerate}
\end{enumerate}

Next, we prove (iii).
 When $p\leq j$, we have
$j+\frac{p^2}{2}-\frac{3p}{2}+1-j=\frac{p^2-3p+2}{2}\geq 0$ for all $p\geq 2$. When $p\geq j+1$,
\begin{align*}
    \frac{j}{2}-\frac{j^2}{2}+pj-p+1-p&=\frac{j}{2}-\frac{j^2}{2}+(j-2)p+1\\
                                      &\geq\frac{j}{2}-\frac{j^2}{2}+(j-2)(j+1)+1\\
                                      &= \textstyle\frac 12(j^2-j-2)\geq 0
    \end{align*}
for all $j\geq 2$. The proof is complete.\qed

Throughout the paper, $U$ shall be used to denote the set Alice marks and $I\subseteq U$ denotes
the set Bob removes. For any $U\subseteq V(G)$, the indicator function $\mathbf{1}_U$ of $U$ is
defined as $\mathbf{1}_U(x)=\begin{cases}1 & \mbox{ if }x\in U;\\ 0 & \mbox{ if } x\not\in
U.\end{cases}$

\begin{proposition} \cite{KMZ13+,Schauz09} If $G$ is edgeless and $f(v)\geq 1$ for all $v\in V(G)$, then $G$ is
on-line $f$-choosable. A graph $G$ is on-line $f$-choosable if and only if for any $U\subseteq
V(G)$ there exists an independent set $I\subseteq U$ of $G$ such that $G-I$ is on-line
$(f-\mathbf{1}_U)$-choosable. \label{pro4}\end{proposition}

For a subset $Y\subseteq X\in \overline{X}_j$ with $|Y|=y$, we say that $F(y)$ is {\it saturated with
respect to $Y$} if $F(y)=\sum_{w\in Y}f(w)=\alpha(j)+\beta(j)$. We are now ready to prove the main result.

\begin{theorem}\label{thm1}
Let $G$ be a complete multipartite graph with independence number $m\geq 2$. If there is a
partition $\Pi$ of parts of $G$ and a function $f: V(G)\rightarrow \mathbb{N}$
satisfying the following:
\begin{enumerate}
\item[{(R1)}] $F(j)\geq \alpha(j)+\beta(j)$ for all $j=1,\ldots,m$ and
\item[{(R2)}] $\displaystyle f(v)\geq S(m-j)+\sum_{p=1}^{j-1}|V(\underline{X}_p)|+\sum_{q=1}^{i-1}|V(\underline{X}_{j,q})|+1$ for all
$v\in V(\underline{X}_{j,i})$ for all $j = 1, ... ,m$ and $1 \leq i \leq \ell_j$,
\end{enumerate}
then $G$ is on-line $f$-choosable.
\end{theorem}
\proof We shall prove the theorem by induction on $|V(G)|$. Obviously, if $G$ is edgeless, then $G$
is on-line $f$-choosable since $f(v)\geq F(1)\geq 1$ for all $v\in V(\bigcup_{p=2}^m\overline{X}_p)$ by
(R1) and $f(v)\geq 1$ for $v\in V(\bigcup_{p=1}^{m-1}\underline{X}_p)$ by (R2). Assume that $G$ has at
least two parts and that the statement is true for all graphs of order less than $|V(G)|$. We shall
prove that if $G$ has a partition $\Pi$ of parts and a function $f: V(G)\rightarrow \mathbb{N}$ so
that (R1) and  (R2) are satisfied, then no matter what $U\subseteq V(G)$ Alice marks, there exists
an independent set $I\subseteq U$ of $G$ such that the resulting graph $G'=G-I$ satisfies the two
conditions, i.e., there exists a partition $\Pi'$ of parts of $G'$ such that $f'=f-\mathbf{1}_U$ satisfies
(R1) and  (R2) with respect to $\Pi'$. Then by induction we conclude that $G'$ is on-line
$f'$-choosable and thus $G$ is on-line $f$-choosable by Proposition \ref{pro4}.

For a given $U\subseteq V(G)$, the crucial step is twofold: decide an independent set $I\subseteq
U$ and give a partition $\Pi'=\{\underline{X}'_{m-1}$, $\underline{X}'_{m-2}, \ldots, \underline{X}'_{1}$, $\overline{X}'_{2},
\ldots, \overline{X}'_{m}\}$ of parts of $G'$. Our strategy will be given case by case
depending on $U$. Particularly, in any considered case we shall assume that all the previous cases
do not hold. Note that from $\Pi$ to $\Pi'$ all families are inherited except two: the family from
which $I$ is chosen and the family where the remaining partite set $X- I$ is inserted.  The
notations $(\underline{X}_{j}, X- I)$ and $(X- I, \underline{X}_{j})$ denote that the remaining set $X- I$ is
inserted to the end and the beginning of the family $\underline{X}_{j}$, respectively.  Note also that,
once $U$ is given, the function $f'$ can be obtained from $f$ with little difference $\mathbf{1}_U$. So we
may and shall verify the inequalities in {\bf (R1)} and {\bf (R2)} for $f'$ and $\Pi'$ by comparing
the difference with that for $f$ and $\Pi$.



{\bf Case 1: $U$ contains a part $X\in \overline{X}_{j^*}$ for some $j^*$.}\\ 
Let $I=X$ and $\Pi'$ be obtained from $\Pi$ by removing $I$ where all families remain except
$\overline{X}'_{j^*}=\overline{X}_{j^*}- X$.
Next we verify (R1) and (R2) for the updated $\Pi'$ and $f'$. \\
{\bf (R1)}. For all $j\in [m]$, we have $F'(j)\geq F(j)-j\geq \alpha(j)+\beta(j)-j\geq
\alpha'(j)+\beta'(j)$, where the last inequality follows from $u'_{j^*}=u_{j^*}-1$, $\beta'(j)=\beta(j)$ and Proposition \ref{pro2}.\\
{\bf (R2)}.  Since $f'(v)\geq f(v)-1$ for each $v\in V(\underline{X}_j)$, it suffices to show that
$S'(m-j)\leq S(m-j)-1$. This follows immediately from $u'_{j^*}=u_{j^*}-1$ and Proposition
\ref{pro1}.

{\bf Case 2:  $U\cap X\neq \emptyset$ for some $X\in \overline{X}_{j^*}$ and $F(y)$ is saturated with
respect to $Y$ for some $Y\subseteq U\cap X$, where $y=|Y|$.}\\ Among all those cases we choose the
one with the largest $y$. Let $I=U\cap X$ and $\Pi'$ be obtained from $\Pi$ by removing $I$ where
all families remain except $\overline{X}'_{j^*}=\overline{X}_{j^*}- X$ and $\underline{X}'_{m-y}=(\underline{X}_{m-y},X- I)$.
Let $y^*=|U\cap X|$. Obviously, $u'_{j^*}=u_{j^*}-1$, $\ell'_{m-y}=\ell_{m-y}+1$ and
$|V(\underline{X}'_{m-y})|=|V(\underline{X}_{m-y})|+(j^*-y^*)$.  Of particular note is that $j^*>y$ otherwise it
is Case 1.
\\
{\bf (R1)}. Consider $F'(j)$ for the case $j\leq y$, which implies $j^*>j$. Since $j^*>j$ and
$u'_{j^*}=u_{j^*}-1$, from Eq.(\ref{eq1}) we obtain
\begin{equation}\label{eq3}\alpha(j)=\alpha'(j)+\left(\frac{j}{2}-\frac{j^2}{2}+jj^*-j^*+1\right).\end{equation} Since
there are $j^*-y^*$ elements inserted to the family $\underline{X}'_{m-y}$ and $m-y\leq m-j$, by
Proposition \ref{pro2}, $\beta(j)=\beta'(j)-(j-1)(j^*-y^*)$. We now verify (R1),
 \begin{align*}
    \d F'(j)&\geq F(j)-j\\
            &\geq \alpha(j)+\beta(j)-j \\
            &=\left[\alpha'(j)+\left(\frac{j}{2}-\frac{j^2}{2}+jj^*-j^*+1\right)\right]+\left[\beta'(j)-(j-1)(j^*-y^*)\right]-j\\
            &= \alpha'(j)+\beta'(j)+(j-1)\left(y^*-\frac{j+2}{2}\right)\\
            &\geq \alpha'(j)+\beta'(j) \hspace{1cm} (\mbox{as } y^*\geq y\geq j).
\end{align*}
Consider the case $j>y$. Since $j^*>y$, it follows that $m-y\geq m-j^*+1$ and the coefficient of
$|V(\underline{X}_{m-y})|$ in the expression (\ref{eq2}) of $\beta(j)$ is $y$. Thus,
$\beta(j)=\beta'(j)-y(j^*-y^*)$. We have $F'(j)\geq F(j)-j \geq
\underbrace{\alpha(j)+\beta(j)-j+1}_{Q}$
where the second inequality holds by the maximality assumption of $y$.\\
If $j<j^*$, then
\begin{align*}Q&=\underbrace{\alpha'(j)+\left(\frac{j}{2}-\frac{j^2}{2}+jj^*-j^*+1\right)}_{\mbox{by Eq.
(\ref{eq3})}}+\beta'(j)-y(j^*-y^*)-j
+1\\ &\geq \alpha'(j)+\left(\frac{j}{2}-\frac{j^2}{2}+jj^*-j^*+1\right)+\beta'(j)-y^*(j^*-y^*)-j +1 \\
&\geq \alpha'(j)+\beta'(j),
\end{align*} where the last inequality can be proved through two cases:
If $y^*\geq j$, then it follows from the same analysis as the previous case; if $y^*<j$, then it
follows from the fact that the quadratic function
$g(j)=\left(\frac{j}{2}-\frac{j^2}{2}+jj^*-j^*+1\right)-y^*(j^*-y^*)-j+1$ is increasing when $y^*<j<j^*$ and
$g(y^*+1)\geq 0$ for $y^*\notin(1,2)$.\\
If $j\geq j^*$, then the coefficient of $u_{j^*}$ in the expression (\ref{eq1}) of $\alpha(j)$ is
$j+\frac{(j^*)^2}{2}-\frac{3j^*}{2}+1$, which implies
$\alpha(j)=\alpha'(j)+\left(j+\frac{(j^*)^2}{2}-\frac{3j^*}{2}+1\right)$. Hence
\begin{align*}Q&=\left[\alpha'(j)+\left(j+\frac{(j^*)^2}{2}-\frac{3j^*}{2}+1\right)\right]+\beta'(j)-y(j^*-y^*)-j
+1\\ &\geq \alpha'(j)+\beta'(j)+\frac{(j^*)^2}{2}-\frac{3j^*}{2}+2-y(j^*-y)\\
&\geq \alpha'(j)+\beta'(j)+\frac{(j^*)^2}{2}-\frac{3j^*}{2}+2-\frac{(j^*)^2}{4} \hspace{0.5cm} \left(\text{as $y(j^*-y)$ has a max. at $y=\frac{j^*}2$}\right)  \\
&= \alpha'(j)+\beta'(j)+\frac{(j^*-2)(j^*-4)}{4}\\
&\geq \alpha'(j)+\beta'(j) \mbox{ whenever } j^*\neq 3 \mbox{ (noticing that } j^*\geq 2).
\end{align*}
When $j^*=3$, either $y=1$ or $y=2$. This implies $y(j^*-y)=2$.
Consequently, $\frac{(j^*)^2}{2}-\frac{3j^*}{2}+2-y(j^*-y)=0$ and then $F'(j)\geq
Q\geq \alpha'(j)+\beta'(j)$, as desired.\\
{\bf (R2)}. Consider $f'(v)$ for any $v\in \underline{X}'_{j,i}$.\\
For the case $\{j<m-y\}$ or $\{j=m-y$ and $i\leq \ell_{m-y}\}$, we conclude the same
\begin{align*}f'(v)\geq
    f(v)-1&\geq S(m-j)+\sum_{p=1}^{j-1}|V(\underline{X}_p)|+\sum_{q=1}^{i-1}|V(\underline{X}_{j,q})|\\
          &\geq S'(m-j)+1+\sum_{p=1}^{j-1}|V(\underline{X}'_{p})|+\sum_{q=1}^{i-1}|V(\underline{X}'_{j,q})|,
\end{align*}
where the last inequality follows from $u'_{j^*}=u_{j^*}-1$, Proposition \ref{pro1} and the fact
that the last two terms are inherited.

For the case $j=m-y$ and $i=\ell_{m-y}+1$, obviously $v\in X- I$. In this case,
\begin{align*}
    f'(v)=f(v)&\geq F(y+1)-\sum_{w\in Y}f(w) \hspace{1cm} \text{(as $F(y)$ is saturated with respect to $Y$)}\\
              &=F(y+1)-[\alpha(y)+\beta(y)]\\ &\geq \alpha(y+1)-\alpha(y)+\beta(y+1)-\beta(y)\\
              &=S(y)+\sum_{p=1}^{m-y}|V(\underline{X}_p)|\\
              &\geq\left[S'(y)+1\right]+\left[\sum_{p=1}^{m-y-1}|V(\underline{X}'_{p})|+\sum_{q=1}^{i-1}|V(\underline{X}_{m-y,q})|\right] \text{(by Proposition \ref{pro1})}  
\end{align*}
and thus (R2) is satisfied by simply substituting $m-j$ into $y$.

For the case $j>m-y$, we have $j^*>y^*\geq y\geq m-j+1$. Therefore, by definition, the coefficient
of $u_{j^*}$ in the expression of $S(m-j)$ is $j^*-m+j$. This implies $S'(m-j)\leq
S(m-j)-(j^*-m+j)$ and thus
\begin{align*}f'(v)&\geq f(v)-1\\ &\geq S(m-j)+\left[\sum_{p=1}^{j-1}|V(\underline{X}_p)|+\sum_{q=1}^{i-1}|V(\underline{X}_{j,q})|\right]\\ &\geq
S'(m-j)+(j^*-m+j)+\left[\sum_{p=1}^{j-1}|V(\underline{X}'_{p})|+\sum_{q=1}^{i-1}|V(\underline{X}'_{j,q})|-(j^*-y^*)\right]\\
&\geq S'(m-j)+\sum_{p=1}^{j-1}|V(\underline{X}'_{p})|+\sum_{q=1}^{i-1}|V(\underline{X}'_{j,q})|+1,\end{align*} as
desired.   

Note that as Cases 1 and 2 were excluded, in all remaining cases we conclude that if $U\cap X\neq
\emptyset$ for some $X\in \overline{X}_j$, then
 \begin{itemize}\item[(i)] $|U\cap X|<j$ otherwise it is Case 1,
 \item[(ii)] for all subsets $Y\subseteq U\cap X$, $F(|Y|)$ is not saturated with
 respect to $Y$ otherwise it is Case 2.
\end{itemize}

{\bf Case 3: From the above discussion, we know that one of the following cases must happen if $U\neq
\emptyset$. \begin{itemize}\item[(3.1)] $V(\underline{X}_{m+1-s})\cap U\neq \emptyset$ for some $2\leq
s\leq m$,
\item[(3.2)] there exists $Y\subseteq U$ and $\displaystyle Y\subseteq X\in \bigcup_{p=2}^m\overline{X}_p$ with $|Y|=t$ for some $1\leq t\leq m-1$ such that  $F(t)$
 is  not saturated with respect to $Y$,
\item[(3.3)] there exists $Y\nsubseteq U$ and $Y\subseteq X\in
\overline{X}_u$ for some $2\leq u\leq m$ such that $F(|Y|)$ is saturated with respect to $Y$ and $X\cap
U\neq \emptyset$.
\end{itemize}}
\noindent Among all these cases, we choose the one with the largest $s,t,u$ (if they are equal then the
priority is $s$, $t$ and then $u$). For example, if $s=u=4$ and $t=5$ then the first case we deal
with is Case (3.2) with $t=5$.

{\bf Case (3.1): Notice that in this case, $s\geq u$ and $s\geq t$.} We choose the least $i^*$
among all $i$'s with $V(\underline{X}_{m+1-s,i})\cap U\neq \emptyset$. Let $I=V(\underline{X}_{m+1-s,{i^*}})\cap
U$ and $\Pi'$ be obtained from $\Pi$ by removing $I$ where all families remain except
$\underline{X}'_{m+1-s,{i^*}}=\underline{X}_{m+1-s,{i^*}}- I$ and the orderings are inherited. Obviously,
$\alpha'(j)=\alpha(j)$ for all $j\in [m]$.
\\
{\bf (R1)}. Consider any subset $J\subseteq X\in \overline{X}_{j^*}$ with $|J|=j$ for some $j^*$. Notice that
$j\leq j^*$. Here we may assume that $J\cap U\neq\emptyset$ because it is trivial that $\sum_{w\in
J}f'(w)\geq \alpha'(j)+\beta'(j)$ if $J\cap U=\emptyset$.

 If $F(j)$ is saturated with respect to $J$, then $J\nsubseteq U$ (by Case $2$) and
$j^*\leq s$ (otherwise it is Case (3.3) with $u=j^*>s$). Since $\underline{X}'_{m+1-s}=\underline{X}_{m+1-s}-
I$, $I\neq \emptyset$, and the coefficient of $|V(\underline{X}_{m+1-s})|$ in the expression of $\beta(j)$
is $j-1$ (from $j\leq s$ and Proposition \ref{pro2}), we have $\beta(j)\geq \beta'(j)+(j-1)$. It
follows that
\begin{align*}
    \sum_{w\in J}f'(w)&\geq F(j)-(j-1) \hspace{1.5cm} \text{(as $J\nsubseteq U$)} \\
                      &\geq \alpha(j)+\beta(j)-(j-1)\\
                      &\geq \alpha'(j)+\beta'(j). \hspace{1cm} \text{(as $\alpha'(j)=\alpha(j)$)}
\end{align*}

If $F(j)$ is not saturated with respect to $J$ and $J\subseteq U$, then $j\leq s$ (otherwise it
is Case (3.2) with $t=j>s$). Thus, the coefficient of $|V(\underline{X}_{m+1-s})|$ in $\beta(j)$ is $j-1$.
This implies that $\beta(j)\geq \beta'(j)+(j-1)$. Accordingly, $\sum_{w\in J}f'(w)\geq
[\alpha(j)+\beta(j)+1]-j\geq \alpha'(j)+\beta'(j)$, where the first inequality holds as $F(j)$ is
not saturated with respect to $J$.

 If $F(j)$ is not saturated and $J\nsubseteq U$, then $F(|J\cap
U|)$ must also not be saturated with respect to $J\cap U$ (by Case 2) and $|J\cap U|\leq s$ (otherwise it is Case (3.2) with $t=|J\cap U|>s$). Hence, $\sum_{w\in J}f'(w)\geq F(j)+1-|J\cap
U|\geq  \alpha'(j)+\beta'(j)$. Notice that the last inequality holds as $\beta(j)+1-|J\cap
U|\geq\beta'(j)$, which follows from the fact that the coefficient of $|V(\underline{X}_{m+1-s})|$ of
$\beta(j)$ in Eq.(\ref{eq2})  is at least $\min\{j-1,s-1\}$.

From the above discussion, in either case we have $\sum_{w\in J}f'(w)\geq \alpha'(j)+\beta'(j)$. As
$J$ is chosen arbitrarily, we can conclude that $F'(j)\geq \alpha'(j)+\beta'(j)$ for all $j\in
[m]$.
\\
{\bf (R2)}. Because of the maximality of $s$, it suffices to consider $v\in V(\underline{X}'_{j,i})$ for
the two cases $\{j>m+1-s\}$ and $\{j=m+1-s$ and $i>i^*\}$ as $f'(v)$ is inherited otherwise. If
$j>m+1-s$, then $f'(v)\geq f(v)-1\geq
S(m-j)+\sum_{p=1}^{j-1}|V(\underline{X}_p)|+\sum_{q=1}^{i-1}|V(\underline{X}_{j,q})| \geq
S'(m-j)+\left[\sum_{p=1}^{j-1}|V(\underline{X}'_{p})|+1\right]+\sum_{q=1}^{i-1}|V(\underline{X}'_{j,q})|$, as
desired.

For the case that $i>i^*$ and $j=m+1-s$, because of the minimality of $i^*$, we have
$\sum_{q=1}^{i-1}|V(\underline{X}_{j,q})|\geq \sum_{q=1}^{i-1}|V(\underline{X}'_{j,q})|+1$. Similarly, we can
conclude that
 $f'(v)\geq S'(m-j)+\sum_{p=1}^{j-1}|V(\underline{X}'_{p})|+\left[\sum_{q=1}^{i-1}|V(\underline{X}'_{j,q})|+1\right]$, as
 desired.

{\bf Case (3.2): Note that in this case, $t$ is the largest, i.e., $t>s$ and $t\geq u$.} We may
assume that there exists $Y\subseteq U$ and $Y\subseteq X\in \overline{X}_{j^*}$ with $|Y|=t$,
$1\leq t\leq m-1$, such that $F(t)$ is not saturated with respect to $Y$.  Let $I=U\cap X$
(noticing that $j^*>|I|\geq t$) and $\Pi'$ be obtained from $\Pi$ by removing $I$ where all
families remain except that $\overline{X}'_{j^*}=\overline{X}_{j^*}- X$, {\bf(a)
$\overline{X}'_{j^*-|I|}=\overline{X}_{j^*-|I|}\cup \{X- I\}$ if $j^*-|I|\geq 2$} and {\bf(b)
$\underline{X}'_{1}=(X- I, \underline{X}_{1})$ if $j^*-|I|= 1$.} In other words, we have
$u'_{j^*}=u_{j^*}-1$ and either (a) $u'_{j^*-|I|}=u_{j^*-|I|}+1$ or (b)
$|V(\underline{X}'_{1})|=|V(\underline{X}_{1})|+1$. Observing the corresponding coefficients of
$u_{j^*}, u_{j^*-|I|}$ and $|V(\underline{X}_{1})|$ in the expression of $\alpha(j)+\beta(j)$, we
have that $a_{j^*}-a_{j^*-|I|}$ and $a_{j^*}-b_1$ are at least $|I|$ by Proposition \ref{pro2}.
Accordingly, we can conclude that $\alpha(j)+\beta(j)\geq \alpha'(j)+\beta'(j)+|I|$.
\\
{\bf (R1)}. Consider any subset $J\subseteq K\in \overline{X}_{k}$ with $|J|=j$ for some $k$. Notice that
$j\leq k$. We may assume that $J\cap U\neq\emptyset$, otherwise $\sum_{w\in J}f'(w)\geq F(j)\geq
\alpha'(j)+\beta'(j)$ since $\alpha(j)+\beta(j)\geq \alpha'(j)+\beta'(j)+|I|$.

For the special case that $K=X- I\in \overline{X}'_{j^*-|I|}$, $\sum_{w\in J}f'(w)=\sum_{w\in J}f(w)$ and
thus (R1) is satisfied immediately as $\alpha'(j)\leq \alpha(j)$ and $\beta'(j)=\beta(j)$.

Consider the case that $F(j)$ is saturated with respect to $J$ and $k\leq t$, which implies
$J\nsubseteq U$ by Case $2$. It follows that $\sum_{w\in J}f'(w)\geq \alpha(j)+\beta(j)-(j-1)$.
 Since $j\leq k\leq t\leq |I|$, $\sum_{w\in J}f'(w)\geq
\alpha(j)+\beta(j)-(j-1)\geq \alpha'(j)+\beta'(j)+|I|-(j-1)> \alpha'(j)+\beta'(j)$.

If $F(j)$ is saturated with respect to $J$ and $k> t$,  then $J\cap U$ must be empty otherwise it
is either $J\subseteq U$ (Case 2) or $J\nsubseteq U$ (Case (3.3) with $u=k>t$), a contradiction to
the assumption $J\cap U\neq\emptyset$. 

If $F(j)$ is not saturated with respect to $J$, then $F(|J\cap U|)$ must not be saturated with
respect to $J\cap U$ (by Case 2) and $|J\cap U|\leq t$ (by the maximality assumption of $t$). It
follows that $\sum_{w\in J}f'(w)\geq [\alpha(j)+\beta(j)+1]-|J\cap U|\geq
\alpha'(j)+\beta'(j)+|I|+1-|J\cap U|\geq \alpha'(j)+\beta'(j)$ as
 $|I|\geq t\geq |J\cap U|$.

From the above discussion, in either case we have $\sum_{w\in J}f'(w)\geq \alpha'(j)+\beta'(j)$. As
$J$ is chosen arbitrarily, we can conclude that $F'(j)\geq \alpha'(j)+\beta'(j)$ for all $j\in
[m]$.
 \\
{\bf (R2)}. Consider $f'(v)$ for any $v\in \underline{X}_{j,i}$.

For the special case that $v\in X- I=\underline{X}_{1,1}$, $f'(v)=f(v)\geq F(1)\geq
\alpha(1)+\beta(1)=\sum_{p=2}^mu_p=S'(m-1)+1$, as desired. Recall that if it is the case (a), all families remain except
that $\overline{X}'_{j^*}=\overline{X}_{j^*}- X$ and $\overline{X}'_{j^*-|I|}=\overline{X}_{j^*-|I|}\cup \{X- I\}.$ If it is the case (b), all families remain except
that $\overline{X}'_{j^*}=\overline{X}_{j^*}- X$ and $\underline{X}'_{1}=(X- I, \underline{X}_{1})$.

 If $j\leq m+1-t$, then we have $v\not\in V(\underline{X}_{j})\cap U$ otherwise it is Case (3.1) with $s=m+1-j\geq t$. Therefore,
$f'(v)=f(v)\geq S(m-j)+\sum_{p=1}^{j-1}|V(\underline{X}_p)|+\sum_{q=1}^{i-1}|V(\underline{X}_{j,q})|+1$.\\
If it is the case (a), then $S'(m-j)\leq S(m-j)$ (as $s_{j^*}\geq s_{j^*-|I|}$ by Proposition
\ref{pro1}) and
$\sum_{p=1}^{j-1}|V(\underline{X}'_{p})|+\sum_{q=1}^{i-1}|V(\underline{X}'_{j,q})|=\sum_{p=1}^{j-1}|V(\underline{X}_p)|+\sum_{q=1}^{i-1}|V(\underline{X}_{j,q})|$.\\
If it is the case (b), then $S'(m-j)\leq S(m-j)-1$ (by Proposition \ref{pro1} and
$u'_{j^*}=u_{j^*}-1$) and
$\sum_{p=1}^{j-1}|V(\underline{X}'_{p})|+\sum_{q=1}^{i-1}|V(\underline{X}'_{j,q})|=\sum_{p=1}^{j-1}|V(\underline{X}_p)|+\sum_{q=1}^{i-1}|V(\underline{X}_{j,q})|+1$.
In either case, $f'(v)\geq
S'(m-j)+\sum_{p=1}^{j-1}|V(\underline{X}'_{p})|+\sum_{q=1}^{i-1}|V(\underline{X}'_{j,q})|+1$, as desired.

If $j> m+1-t$, then $m+1-j<j^*$ since $t<j^*$ otherwise it is Case 1. In this case, observing
the coefficients of $u_{j^*}$ and $u_{j^*-|I|}$ in $S(m-j)$, we obtain $s_{j^*}>s_{j^*-|I|}$ and
$s_{j^*}\geq 2$ by Proposition
\ref{pro1}. \\
If it is the case (a), then $S(m-j)\geq S'(m-j)+1$ (since $s_{j^*}>s_{j^*-|I|}$) and thus
$f'(v)\geq f(v)-1\geq
S'(m-j)+\sum_{p=1}^{j-1}|V(\underline{X}'_{p})|+\sum_{q=1}^{i-1}|V(\underline{X}'_{j,q})|+1$
as $\sum_{p=1}^{j-1}|V(\underline{X}'_{p})|$ and $\sum_{q=1}^{i-1}|V(\underline{X}'_{j,q})|$ are inherited.\\
If it is the case (b), then $S(m-j)\geq S'(m-j)+2$ (since $s_{j^*}\geq 2$ and $u'_{j^*}=u_{j^*}-1$)
and
$\sum_{p=1}^{j-1}|V(\underline{X}'_{p})|+\sum_{q=1}^{i-1}|V(\underline{X}'_{j,q})|=\sum_{p=1}^{j-1}|V(\underline{X}_p)|+\sum_{q=1}^{i-1}|V(\underline{X}_{j,q})|+1$.
Thus, $f'(v)\geq f(v)-1\geq
[S'(m-j)+2]+\left[\sum_{p=1}^{j-1}|V(\underline{X}'_{p})|+\sum_{q=1}^{i-1}|V(\underline{X}'_{j,q})|-1\right]\geq
S'(m-j)+\sum_{p=1}^{j-1}|V(\underline{X}'_{p})|+\sum_{q=1}^{i-1}|V(\underline{X}'_{j,q})|+1$, as desired.

{\bf Case (3.3): Note that in this case, $u$ is the largest, i.e., $u>s$ and $u>t$.} To avoid
confusion, instead of $u$ we shall use $j^*$ to denote the largest index and assume that there
exists $Y\nsubseteq U$ and $Y\subseteq X\in \overline{X}_{j^*}$ such that $F(|Y|)$ is saturated with
respect to $Y$.  Among all these cases with the same $j^*$, we choose the one with the largest
$|U\cap X|$. Let $I=U\cap X$ and $\Pi'$ be obtained from $\Pi$ by removing $I$ where all families
remain except that $\overline{X}'_{j^*}=\overline{X}_{j^*}- X$, {\bf (a)
$\overline{X}'_{j^*-|I|}=\overline{X}_{j^*-|I|}\cup \{X- I\}$ if $j^*-|I|\geq 2$} and {\bf (b) $\underline{X}'_{1}=(X-
I,\underline{X}_{1})$ if $j^*-|I|=1$.} Notice that $u'_{j^*}=u_{j^*}-1$ and either (a)
$u'_{j^*-|I|}=u_{j^*-|I|}+1$ or (b) $|V(\underline{X}'_{1})|=|V(\underline{X}_{1})|+1$. The same argument as that
in Case (3.2) implies that $\alpha(j)+\beta(j)\geq \alpha'(j)+\beta'(j)+|I|$.
\\
{\bf (R1)}. Consider any subset $J\subseteq K\in \overline{X}_{k}$ with $|J|=j$ for some $k$  and $J\cap
U\neq\emptyset$. Notice that $j\leq k$.

For the special case that $K=X- I\in \overline{X}'_{j^*-|I|}$, it follows that $\sum_{w\in
J}f'(w)=\sum_{w\in J}f(w)$ and thus (R1) is satisfied immediately as $\alpha'(j)\leq \alpha(j)$ and
$\beta'(j)=\beta(j)$.

If  $F(j)$ is saturated with respect to $J$, then $k\leq j^*$ otherwise it is Case (3.3) with
$k>j^*$ violating the maximality assumption of $j^*$. Next we discuss the two cases $k=j^*$ and
$k<j^*$ separately.\\ For the case $k=j^*$, $\sum_{w\in J}f'(w)\geq F(j)-|I|\geq
\alpha(j)+\beta(j)-|I|\geq \alpha'(j)+\beta'(j)$ where the first inequality follows from the fact
that among all the cases with $k=j^*$ we choose the one with the largest $|U\cap X|$. \\If $k<j^*$,
then $J\nsubseteq U$ (since Case 2 does not hold). It follows that $\sum_{w\in J}f'(w)\geq
\alpha(j)+\beta(j)-(j-1)$. In addition, $j\leq k<j^*$ implies that
$a_{j^*}-\max\{a_{j^*-|I|},b_1\}\geq j-1$
 by Proposition \ref{pro2}. It follows that $\alpha(j)+\beta(j)\geq
\alpha'(j)+\beta'(j)+j-1$. Hence, $\sum_{w\in J}f'(w)\geq \alpha'(j)+\beta'(j)$.

Consider the case when $F(j)$ is not saturated with respect to $J$. It follows that $\sum_{w\in J}f'(w)\geq \alpha(j)+\beta(j)+1-j$.
There are two cases: $j< j^*$ and $j\geq j^*$.\\
If $j<j^*$, by Proposition \ref{pro2} we have $a_{j^*}-\max\{a_{j^*-|I|},b_1\}\geq j-1$. This
implies that $\alpha(j)+\beta(j)\geq \alpha'(j)+\beta'(j)+j-1$ and therefore $\sum_{w\in
J}f'(w)\geq \alpha'(j)+\beta'(j)$.\\
If  $j\geq j^*$, then we have $J\nsubseteq U$ otherwise it is Case (3.2) with $t=j\geq j^*=u$.
Furthermore, it must be $|J\cap U|<j^*$ otherwise either $F(|J\cap U|)$ is saturated with
respect to $J\cap U$ (Case 2) or $F(|J\cap U|)$ is not saturated with respect to $J\cap U$ and
$t=|J\cap U|\geq j^*=u$ (Case (3.2)). Since $F(j)$ is not saturated with respect to $J$ and $|J\cap
U|<j^*$, we have $\sum_{w\in J}f'(w)\geq \alpha(j)+\beta(j)+1-(j^*-1)$. Since $j\geq j^*$, by
Proposition \ref{pro2} we have $a_{j^*}-\max\{a_{j^*-|I|},b_1\}\geq j^*-2$. It follows that
$\alpha(j)+\beta(j)\geq \alpha'(j)+\beta'(j)+j^*-2$ and thus $\sum_{w\in J}f'(w)\geq
\alpha'(j)+\beta'(j)$.

From the above discussion, in either case we have $\sum_{w\in J}f'(w)\geq \alpha'(j)+\beta'(j)$. As
$J$ is chosen arbitrarily, we can conclude that $F'(j)\geq \alpha'(j)+\beta'(j)$ for all $j\in
[m]$.
\\
{\bf (R2)}. Consider $f'(v)$ for any $v\in \underline{X}_{j,i}$.

For the special case that $v\in X- I=\underline{X}_{1,1}$, $f'(v)=f(v)\geq F(1)\geq
\alpha(1)+\beta(1)=\sum_{p=2}^mu_p=S'(m-1)+1$, as desired. Recall that if it is the case (a), all families remain except
that $\overline{X}'_{j^*}=\overline{X}_{j^*}- X$ and $\overline{X}'_{j^*-|I|}=\overline{X}_{j^*-|I|}\cup \{X- I\}$. If it is the case (b), all families remain except
that $\overline{X}'_{j^*}=\overline{X}_{j^*}- X$ and $\underline{X}'_{1}=(X-
I,\underline{X}_{1})$.

 If $j\leq m+1-j^*$, then we have $v\not\in V(\underline{X}_{j})\cap U$ otherwise it is Case (3.1) with $s=m+1-j\geq j^*=u$. Therefore,
$f'(v)=f(v)\geq S(m-j)+\sum_{p=1}^{j-1}|V(\underline{X}_p)|+\sum_{q=1}^{i-1}|V(\underline{X}_{j,q})|+1$.\\
If it is the case (a), then $S'(m-j)\leq S(m-j)$ (by Proposition \ref{pro1} and $j^*\geq j^*-|I|$)
and
$\sum_{p=1}^{j-1}|V(\underline{X}'_{p})|+\sum_{q=1}^{i-1}|V(\underline{X}'_{j,q})|=\sum_{p=1}^{j-1}|V(\underline{X}_p)|+\sum_{q=1}^{i-1}|V(\underline{X}_{j,q})|$.\\
If it is the case (b), then $S'(m-j)\leq S(m-j)-1$ (by Proposition \ref{pro1} and
$u'_{j^*}=u_{j^*}-1$) and
$\sum_{p=1}^{j-1}|V(\underline{X}'_{p})|+\sum_{q=1}^{i-1}|V(\underline{X}'_{j,q})|=\sum_{p=1}^{j-1}|V(\underline{X}_p)|+\sum_{q=1}^{i-1}|V(\underline{X}_{j,q})|+1$.\\
In either case, $f'(v)\geq
S'(m-j)+\sum_{p=1}^{j-1}|V(\underline{X}'_{p})|+\sum_{q=1}^{i-1}|V(\underline{X}'_{j,q})|+1$ as desired.

If $j> m+1-j^*$, then $j^*>m+1-j$.
If it is the case (a), then to verify $f'(v)$ it suffices to compare the coefficients of $u_{j^*}$
and $u_{j^*-|I|}$ in $S(m-j)$. Since $j^*>m+1-j$, we have $s_{j^*}>s_{j^*-|I|}$ by Proposition
\ref{pro1}. Consequently, $S'(m-j)\leq S(m-j)-1$. Thus,
$f'(v)\geq S'(m-j)+\sum_{p=1}^{j-1}|V(\underline{X}'_{p})|+\sum_{q=1}^{i-1}|V(\underline{X}'_{j,q})|+1$.\\
Otherwise, if it is the case (b), then
$\sum_{p=1}^{j-1}|V(\underline{X}'_{p})|+\sum_{q=1}^{i-1}|V(\underline{X}'_{j,q})|=\sum_{p=1}^{j-1}|V(\underline{X}_p)|
+\sum_{q=1}^{i-1}|V(\underline{X}_{j,q})|+1$. Since $j^*>m+1-j$, the term $u_{j^*}$ in $S(m-j)$ has a
coefficient $s_{j^*}\geq 2$ by Proposition \ref{pro1}. It follows that $S(m-j)\geq S'(m-j)+2$ as
$u'_{j^*}=u_{j^*}-1$. Thus, $f'(v)\geq f(v)-1\geq
S(m-j)+\sum_{p=1}^{j-1}|V(\underline{X}_p)|+\sum_{q=1}^{i-1}|V(\underline{X}_{j,q})|\geq
(S'(m-j)+2)+(\sum_{p=1}^{j-1}|V(\underline{X}'_{p})|+\sum_{q=1}^{i-1}|V(\underline{X}'_{j,q})|-1)\geq
S'(m-j)+\sum_{p=1}^{j-1}|V(\underline{X}'_{p})|+\sum_{q=1}^{i-1}|V(\underline{X}'_{j,q})|+1.$

As one of the above cases  must occur, by induction the proof is complete. \qed

%
%
%
%


\section{Consequences of Theorem \ref{thm1}}

Theorem \ref{thm1} provides a sufficient condition on $f$ for graphs being on-line $f$-choosable.
It is a widely applicable tool for computing the on-line choice number of complete multipartite
graphs with varying parameters. Many interesting results can be obtained immediately from Theorem
\ref{thm1}. This section provides only some of them that are relevant to recent results.

\begin{theorem}\label{thm2}
Let $G$ be a complete multipartite graph with independence number $m\geq 2$ and let $k_p$ denote the number of parts of cardinality $p$ for $1\leq p\leq m$.
If $\displaystyle k_1-\sum_{p=2}^m\left(\frac{p^2}{2}-\frac{3p}{2}+1\right)k_p\geq 0$, then $G$ is on-line chromatic-choosable.
\end{theorem}

\proof Consider the particular partition $\Pi=\{\underline{X}_{m-1}$, $\underline{X}_{m-2}, \ldots, \underline{X}_{1}$,
$\overline{X}_{2}, \ldots, \overline{X}_{m}\}$ of parts of $G$ where $\underline{X}_{m-1}$, $\underline{X}_{m-2}, \ldots, \underline{X}_{2}$ are
empty sets. Note that when all the families $\underline{X}_{m-1}$, $\underline{X}_{m-2}, \ldots, \underline{X}_{2}$ are
empty, the remaining families are determined exactly. That is $\ell_{m-1}=\ell_{m-2}=\cdots
=\ell_{2}=0$, $\ell_1=k_1$ and $u_p=k_p$ for $p=2,3,\ldots,m.$ Obviously, $\sum_{p=2}^mu_p+\ell_1=\sum_{p=1}^mk_p=\chi(G)$. Let $f(v)=\chi(G)$ for all $v\in V(G)$.

Next, we will verify that $f$ and $\Pi$ with arbitrarily ordered $\underline{X}_{1}$ satisfy (R1) and (R2)
in Theorem \ref{thm1}. The inequality (R2) holds as $f(v)=\sum_{p=2}^mu_p+\ell_1 \geq
S(m-1)+\sum_{q=1}^{i-1}|V(\underline{X}_{1,q})|+1$ for all $v\in V(\underline{X}_{1,i})$.

Consider any subset $J\subseteq X\in \bigcup_{p=2}^m\overline{X}_{p}$ with $|J|=j$, $1\leq j\leq m$. To
verify (R1), it suffices to prove that $\sum_{w\in J}f(w)\geq \alpha(j)+\beta(j)$, or equivalently
that
\begin{multline}
j\left(\sum_{p=2}^mu_p+\ell_1\right)\geq\sum_{p=2}^j\left(j+\frac{p^2}{2}-\frac{3p}{2}+1\right)u_p\\
                                        +\sum_{p=j+1}^m\left(\frac{j}{2}-\frac{j^2}{2}+pj-p+1\right)u_p+(j-1)\ell_1\\
\Leftrightarrow \hspace{.5cm} \ell_1-\sum_{p=2}^j\left(\frac{p^2}{2}-\frac{3p}{2}+1\right)u_p-\sum_{p=j+1}^m\left(-\frac{j}{2}-\frac{j^2}{2}+pj-p+1\right)u_p\geq 0.
\label{eq4'}
\end{multline}
Obviously, Eq.(\ref{eq4'}) is always true for $j=1$. By using elementary calculus, it is easy to
prove that for all $j=2,3,\ldots, m$ with $j\leq p-1$,
$(\frac{p^2}{2}-\frac{3p}{2}+1)-(-\frac{j}{2}-\frac{j^2}{2}+pj-p+1)\geq 0.$ Thus, to prove
Eq.(\ref{eq4'}) for all $j\in [m]$, it suffices to show that
\begin{equation}\ell_1-\sum_{p=2}^m\left(\frac{p^2}{2}-\frac{3p}{2}+1\right)u_p\geq 0.\label{eq5'}\end{equation} This is trivially
true as $k_1-\sum_{p=2}^m\left(\frac{p^2}{2}-\frac{3p}{2}+1\right)k_p\geq 0$. By Theorem
\ref{thm1}, $G$ is on-line chromatic-choosable. \qed

For any two graphs $G$ and $H$, denote by $G+ H$ the {\it join} of $G$ and $H$, that is, the
disjoint union of $G$ and $H$ with the edges $\{uv: u\in V(G), v\in V(H)\}$. Kozik, Micek and Zhu
\cite{KMZ13+} proved that for any graph $G$, the join of $G$ and a complete graph of
 order $|V(G)|^2$ is on-line chromatic-choosable. Later, Carraher, Loeb, Mahoney, Puleo, Tsai and
 West \cite{West13+} improved upon $|V(G)|^2$ with an additional assumption. Precisely, they proved that
for every $d$-degenerate graph $G$ having an optimal proper coloring with color classes of size at
most $m$, $G+ K_t$ is on-line chromatic-choosable if $t\geq (m+1)d$. Theorem \ref{thm2} also
provides an alternative result in this aspect.

\begin{corollary}
Let $G$ be a graph of independence number at most $m$ and have an optimal proper coloring where
there are $k_i$ color classes of cardinality $i$. If \[t\geq
\sum_{p=2}^m\left(\frac{p^2}{2}-\frac{3p}{2}+1\right)k_p-k_1,\] then $G+ K_t$ is on-line
chromatic-choosable.
\end{corollary}

Kozik, Micek and Zhu \cite{KMZ13+} commented that when $|V(G)|\leq \chi(G)+\sqrt{\chi(G)}$, then $G$ is
on-line chromatic-choosable. Later, Carraher et al.  \cite{West13+} showed that the same conclusion
holds under a relaxed condition $|V(G)|\leq \chi(G)+2\sqrt{\chi(G)-1}$. They proposed a weak
version of Conjecture \ref{online}:
\begin{conjecture}[Weak On-Line Ohba's Conjecture]\cite{West13+}
There is a constant $c\in (1,2]$ such that $\chi_p(G)=\chi(G)$ whenever $|V(G)|\leq
c\chi(G)$.\end{conjecture} The weak on-line Ohba's conjecture is still open, to the best of our
knowledge. Following the same argument in Theorem \ref{thm2}, we obtain the following result, which
goes one further step towards the weak conjecture.

\begin{corollary}\label{cor}
If $G$ is a graph with independence number $m\geq 2$ and \[|V(G)|\leq
\frac{m^2-m+2}{m^2-3m+4}\chi(G),\] then $\chi(G)=\chi_p(G)$.
\end{corollary}

\proof As $\displaystyle|V(G)|\leq \frac{m^2-m+2}{m^2-3m+4}\chi(G)$, $\chi(G)=\ell_1+\sum^m_{p=2}u_p$, and $|V(G)|=\ell_1+\sum^m_{p=2}pu_p$, we have $$\displaystyle \ell_1\geq \sum_{p=2}^m
\left[\frac{p(m^2-3m+4)-(m^2-m+2)}{2m-2}\right]u_p.$$ According to Eq.(\ref{eq5'}), it suffices to
show that when  $2\leq p\leq m$, \begin{equation}\displaystyle\frac{p(m^2-3m+4)-(m^2-m+2)}{2m-2}-
\left(\frac{p^2}{2}-\frac{3p}{2}+1\right)\geq 0.\label{eq5}\end{equation} Simplifying yields
$(m-1)p^2-(m^2+1)p+m^2+m\leq 0$. Consider the quadratic function $g(p)=(m-1)p^2-(m^2+1)p+m^2+m$.
Obviously, $g(2)\leq 0$ and $g(m)\leq 0$ whenever $m\geq 2$. Thus, one can conclude that $g(p)\leq
0$ whenever $2\leq p\leq m$. This completes the proof.\qed

Consequently, any graph $G$ with $|V(G)|\leq 2\chi(G)$ and independence number $m\leq 3$ is on-line
chromatic-choosable; the same conclusion was proved independently in \cite{HWZ11,KKLZ12} for $m=2$
and in  \cite{KMZ13+} for $m\leq 3$.

Alon \cite{alon92} established the asymptotically tight bound $\chi_{\ell}(K_{m\star
k})=\Theta(k\log m)$. The following result, which is another immediate consequence of Theorem
\ref{thm1}, gives a general upper bound for $\chi_p(K_{m\star k})$.

\begin{corollary}\label{cor3}
For any integer $m\geq 3$, $\chi_p(K_{m\star k})\leq \left(m+\frac{1}{2}-\sqrt{2m-2}\right)k$.
\end{corollary}

\proof Consider the partition $\Pi=\{\underline{X}_{m-1}$, $\underline{X}_{m-2}, \ldots, \underline{X}_{1}$,
$\overline{X}_{2}, \ldots, \overline{X}_{m}\}$ of parts of $K_{m\star k}$ where all families are empty except
$\overline{X}_{m}$, i.e., $u_m=k$ and $\ell_1=\cdots=\ell_{m-1}=u_2=\cdots=u_{m-1}=0$. Let
$f(v)=\left(m+\frac{1}{2}-\sqrt{2m-2}\right)k$ for all $v\in V(K_{m\star k})$. The corollary
follows directly from Theorem \ref{thm1} with the specified $(\Pi, f)$. \qed

Particularly, when
$m=3$, it assures the same conclusion in \cite{KMZ13+} that $\chi_p(K_{3\star k})\leq
\frac{3}{2}k$.

%
%
%
%

%
%
%

\subsection*{Acknowledgements}
\medskip

The authors would like to thank the anonymous referees for their valuable
suggestions.


\end{document}